# ASYMPTOTICS OF THE MODULE OF A DEGENERATING CONDENSER AND SOME OF THEIR APPLICATIONS

**V. N. Dubinin**

*The well-known asymptotic formula for the module of a condenser with one of the plates degenerating to a point is generalized to the case of a condenser of general type. The condensers under consideration consist of $n$ plates, $n \geq 2$, and the potential functions of condensers take values of different signs on the plates. The asymptotics are considered when one of the plates is fixed while the other $n-1$ plates are constracted to points. Applications of the formula to geometric function theory are given. Among them are inequalities for complex numbers and Green functions and also theorems on the extremal decomposition and distortion theorems for univalent functions. Bibliography: 11 titles.*

**Dedicated to the 90th anniversary of G. M. Goluzin's birth**

Let $\overline{\mathbf{C}}_z$ be the complex sphere. Let $B$ be a domain of $\mathbf{C}_z$ having a Green function, and let $z_0$ be a finite point of $B$. Consider the condenser

$$C(r, B) = (\overline{\mathbf{C}}_z \setminus B, \{z : |z - z_0| \leq r\},$$

where $r > 0$ is sufficiently small. The second plate of the condenser degenerates to a point as $r \to 0$, and we have the following asymptotic formula for the module of the condenser $C(r, B)$:

$$|C(r,B)| = -\frac{1}{2\pi} \log r + M(B, z_0) + o(1), \quad r \to 0. \tag{1}$$

The constant $M(B, z_0)$ is called the *reduced module* of the domain $B$ with respect to the point $z_0$. There are numerous applications of formula (1) in geometric function theory that were initiated by Hayman's monograph [1] and are connected first of all with the application of the method of symmetrization. In the present paper, we study the asymptotics of a generalized condenser which consists of several plates, and the potential function of this condenser can take values of different signs. Thus, we continue our studies begun in [2, 3]. In §1, we obtain an analog of formula (1) in a very general case. In §2, we give some applications connected with the use of a separating transformation and dissymmetrization [4]. We mainly restrict our consideration to the problems in which the extremal configuration possesses a symmetry of the $n$th order. We obtain inequalities for complex numbers, Green functions, and inner radii and some distortion theorems in the theory of univalent functions.

## §1. An asymptotic formula

Everywhere below, $B$ is an open set on $\overline{\mathbf{C}}_z$ whose complement $E_0 = \overline{\mathbf{C}}_z \setminus B$ has positive harmonic measure [5]. Let $z_l$, $l = 1, \dots, m$, be distinct points of the set $B$, let $\delta_l$, $l = 1, \dots, m$, be arbitrary nonzero real numbers, and let $\mu_l$, $\nu_l$, $l = 1, \dots, m$, be arbitrary positive numbers. We use the following notation:

$$Z = \{z_l\}, \quad \Delta = \{\delta_l\}, \quad \Psi = \{\psi_l\}, \quad \psi_l = \psi_l(r) \equiv \mu_l r^{\nu_l}, \quad l = 1, \dots, m;$$

$$E(z_0, r) = \begin{cases} \{z : |z - z_0| \leq r\} & \text{if } z_0 \text{ is a finite point,} \\ \{z : |z| \geq 1/r\} & \text{if } z_0 = \infty; \end{cases} \quad r > 0.$$

Here and below, if not otherwise stipulated, the symbols $\{\ \}$, $\sum$, and $\prod$ respectively denote the collection, summation, and product with respect to all possible indices indicated in the context, with exception of those where the summand in $\sum$ is either equal to $\infty$ or indefinite, and the factor in the product either vanishes or is equal to $\infty$. For sufficiently small $r > 0$, we define a generalized condenser as the ordered collection

$$C(r; B, Z, \Delta, \Psi) = \{E_0, E(z_1, \psi_1(r)), \dots, E(z_m, \psi_m(r))\}$$

with prescribed values 0, $\delta_1, \dots, \delta_m$, respectively. By analogy with the usual condensers, we define the module of the condenser $C(r; B, Z, \Delta, \Psi)$ as the quantity

$$|C(r; B, Z, \Delta, \Psi)| = \sup\left\{1/\iint\limits_{\overline{\mathbf{C}}_z} |\nabla v|^2\, dxdy\right\},$$

where the supremum is taken over all real-valued functions $v(z)$ that are continuous on $\overline{\mathbf{C}}_z$, satisfy the Lipschitz condition in the vicinity of each finite point, vanish in a neighborhood of the set $E_0$, and are such that $v(z) = \delta_l$ on $E(z_l, \psi_l(r))$, $l = 1, \dots, m$ (here $1/0 = +\infty$). The *reduced module* of a set $B$ with respect to the collections $Z$, $\Delta$, and $\Psi$ is defined as the limit*

$$M(B, Z, \Delta, \Psi) = \lim_{r\to 0}\left(|C(r; B, Z, \Delta, \Psi)| + \frac{\nu}{2\pi}\log r\right), \tag{2}$$

where $\nu = (\sum \delta_l^2 \nu_l^{-1})^{-1}$. In the case where $\psi_l(r) \equiv r$, $l = 1, \dots, m$, we write $M(B, Z, \Delta)$ instead of $M(B, Z, \Delta, \Psi)$. Below we show that this limit always exists and is finite, and for now we only note that the definition of the reduced module immediately implies the following important monotonicity property:

$$M(B, Z, \Delta, \Psi) \le M(B', Z, \Delta, \Psi)$$

for any $Z = \{z_l\}$, $z_l \in B \subset B'$, $l = 1, \dots, m$, and any $\Delta$ and $\Psi$. We need auxiliary sets $\widetilde{E}(z_l, r)$, $l = 1, \dots, m$, defined as arbitrary closed subsets of $\overline{\mathbf{C}}_z$ satisfying the condition

$$E(z_l, \psi_l(r_1)) \subset \widetilde{E}(z_l, r) \subset E(z_l, \psi_l(r_2)), \quad l = 1, \dots, m, \tag{3}$$

for some continuous positive $r_j = r_j(r) \sim r$, $r \to 0$, $0 < r < r_0$, $j = 1, 2$. Let $\widetilde{C}(r, B, Z, \Delta, \Psi)$ denote the generalized condenser defined as before, but with $E(z_l, \psi_l(r))$ replaced by $\widetilde{E}(z_l, r)$.

**Lemma 1.** *If for some $\widetilde{E}(z_l, r)$ $l = 1, \dots, m$, there exists a limit*

$$\lim_{r\to 0}\left(|\widetilde{C}(r; B, Z, \Delta, \Psi)| + \frac{\nu}{2\pi}\log r\right), \tag{4}$$

*then there exists a limit (2), and, vice versa, if there exists a limit (2), then there exists a limit (4) equal to it, for any $\widetilde{E}(z_l, r)$, $l = 1, \dots, m$.*

*Proof.* First, we prove that there exists limit (2). In view of continuity, for every sufficiently small $r > 0$ there exists $r'$ and $r''$ such that $r_2(r') = r = r_1(r'')$. Taking (3) into account, we obtain $\widetilde{E}(z_l, r') \subset E(z_l, \psi_l(r)) \subset \widetilde{E}(z_l, r'')$, $l = 1, \dots, m$.

It follows from the definition of the module of a condenser that

$$|\widetilde{C}(r'; B, Z, \Delta, \Psi)| \ge |\widetilde{C}(r; B, Z, \Delta, \Psi)| \ge |\widetilde{C}(r''; B, Z, \Delta, \Psi)|. \tag{5}$$

Note that $r = r_2(r') \sim r'$ and $r = r_1(r'') \sim r''$ as $r \to 0$. Adding $(\nu/2\pi)\log r$ to all parts of inequalities (5) and adding $\pm(\nu/2\pi)\log r'$ (respectively, $\pm(\nu/2\pi)\log r''$) to the left-hand (respectively, right-hand) part, passing to the limit as $r \to 0$, and taking (4) into account, we prove that there exists a limit (2) equal to (4). Vice versa, condition (3) implies that

$$|C(r_1; B, Z, \Delta, \Psi)| \ge |\widetilde{C}(r; B, Z, \Delta, \Psi)| \ge |\, C(r_2; B, Z, \Delta, \Psi)|,$$

where $r_j = r_j(r) \sim r$, $r \to 0$, $j = 1, 2$. Hence, as before, the existence of limit (2) implies the existence of limit (4) equal to it. The lemma is proved.

For an arbitrary domain $C \subset \overline{\mathbf{C}}_z$, let $g_G(z, z_0)$ denote its Green function having a pole at the point $z_0 \in G$ and continued by zero outside of $G$. Let $r(G, z_0)$ be the inner radius of the domain $G$ with respect to the point $z_0$. In the case where $z_0 = \infty$, set $r(G, \infty) = \exp\{\lim\limits_{z\to\infty}(g_G(z, \infty) - \log|z|)\}$.

*The notation used in the present paper differs somewhat from that used in [2].

**Theorem 1.** *Let the set $B$ and the collections $Z$, $\Delta$, and $\Psi$ be the same as before, and let $B_l$ denote the connected component of $B$ containing the point $z_l$, $l = 1, \dots, m$.** *Then*

$$|C(r; B, Z, \Delta, \Psi)| = -\frac{\nu}{2\pi} \log r + M(B, Z, \Delta, \Psi)| + o(1), \quad r \to 0, \tag{6}$$

*where*

$$M(B, Z, \Delta, \Psi) = \frac{\nu^2}{2\pi} \left\{ \sum \frac{\delta_l^2}{\nu_l^2} \log \frac{r(B_l, z_l)}{\mu_l} + \sum \frac{\delta_l \delta_j}{\nu_l \nu_j} g_{B_l}(z_j, z_l) \right\}.$$

*Proof.* First, assume that each of the domains $B_l$ possesses a classical Green function $g_l(z) \equiv g_{B_l}(z, z_l)$, i.e., a Green function which is continuous in $\overline{\mathbf{C}}_z \setminus \{z_0\}$ and vanishes on the complement $\overline{\mathbf{C}}_z \setminus B_l$. Consider an auxiliary function

$$g(z) = \sum_{l=1}^{m} \delta_l \sum_{k=1}^{m} \beta_{kl} g_k(z),$$

where

$$\beta_{kl} = \begin{cases} -(\log \psi_k)^{-1} (\log \psi_l)^{-1} g_k(z_l), & \text{for } k \neq l, \\ -(\log \psi_l)^{-1} [1 + (\log \psi_l)^{-1} \log R_l], & \text{for } k = l, \end{cases}$$

$$R_l = r(B_l, z_l), \quad k = 1, \dots, m; \quad l = 1, \dots, m.$$

The function $g(z)$ vanishes on $\overline{\mathbf{C}}_z \setminus B$ and is harmonic in $B \setminus \{z_l\}$. Set

$$\widetilde{E}(z_l, r) = \{z \in E(z_l, \psi_l(2r)) : g(z)/\delta_l \geq 1\}, \quad l = 1, \dots, m.$$

Since $g(z)$ is harmonic, the boundary of the set $\widetilde{E}(z_l, r)$ with sufficiently small $r$ consists of the points $\{z \in E(z_l, \psi_l(2r)) : g(z) = \delta_l\}$. Indeed, for such points we have for $z_l \neq \infty$ and $r \to 0$

$$\delta_l = \delta_l \left\{ -(\log \psi_l)^{-1} [1 + (\log \psi_l)^{-1} \log R_l](-\log |z - z_l| + \log R_l + o(1)) + o((\log r)^{-1}) \right\}$$

$$+ \sum_{\substack{j=1, \\ j \neq l}}^{m} \delta_j \left\{ -(\log \psi_j)^{-1} [1 + (\log \psi_j)^{-1} \log R_j](g_j(z_l) + o(1)) \right.$$

$$\left. -(\log \psi_l)^{-1} (\log \psi_j)^{-1} g_l(z_j)(-\log |z - z_l| + \log R_l + o(1)) + o((\log r)^{-1}) \right\}.$$

In particular, this implies that $\log |z - z_l| \sim \log \psi_l$ as $r \to 0$. Using the last relation, we obtain

$$\delta_l = \delta_l \left\{ (\log \psi_l)^{-1} \log |z - z_l| - (\log \psi_l)^{-1} \log R_l + (\log \psi_l)^{-2} \log R_l \log |z - z_l| \right\}$$

$$+ \sum_{\substack{j=1, \\ j \neq l}}^{m} \delta_j \left\{ -(\log \psi_j)^{-1} g_j(z_l) + (\log \psi_l)^{-1} (\log \psi_j)^{-1} g_j(z_l) \log |z - z_l| + o((\log r)^{-1}) \right\}$$

$$= \delta_j \left\{ (\log \psi_l)^{-1} \log |z - z_l| + o((\log r)^{-1}) \right\}, \quad r \to 0.$$

Consequently, $|z - z_l| \sim \psi_l$ as $r \to 0$. Hence, the boundary of the set $\widetilde{E}(z_l, r)$ is an "almost circle" with center at the point $z_l$ and radius $\psi_l$, the set $\widetilde{E}(z_l, r)$ itself satisfying condition (3). It is similarly shown that in the case $z_l = \infty$ the boundary of $\widetilde{E}(z_l, r)$ is an "almost circle" with center at the origin and radius $1/\psi_l$, and this set also satisfies condition (3).

*Under such indexing, a component may have several designations.

Now let us consider the condenser

$$\widetilde{C}(r; B, Z, \Delta, \Psi) = \{\overline{\mathbf{C}}_z \setminus B, \widetilde{E}(z_l, r), \dots, \widetilde{E}(z_m, r)\}.$$

By the Dirichlet principle, we have

$$|\widetilde{C}(r; B, Z, \Delta, \Psi)|^{-1} = \iint\limits_{B\setminus \bigcup_{l=1}^{m} \widetilde{E}(z_l,r)} |\nabla g|^2 \, dxdy.$$

Let $\rho > 0$ be so small that the sets $E(z_l, \rho)$ respectively belong to $\widetilde{E}(z_l, r)$, $l = 1, \dots, m$. Applying the Green formula and the Gauss theorem, we subsequently obtain

$$\iint\limits_{B\setminus \bigcup_{l=1}^{m} \widetilde{E}(z_l,r)} |\nabla g|^2 \, dxdy = -\sum_{j=1}^{m} \int\limits_{\partial \widetilde{E}(z_j,r)} \delta_j \frac{\partial g}{\partial n} ds = -\sum_{j=1}^{m} \delta_j \int\limits_{\partial \widetilde{E}(z_j,r)} \frac{\partial g}{\partial n} ds = 2\pi \sum_{j=1}^{m} \delta_j \sum_{l=1}^{m} \delta_l \beta_{jl} + o(1), \ \rho \to 0.$$

Here the differentiation is performed with respect to the inner normal of the field of the condenser $\widetilde{C}(r, B, Z, \Delta, \Psi)$. Substituting the values $\beta_{jl}$ and taking into account that the module of the condenser does not depend on $\rho$, we obtain

$$|\widetilde{C}(r, B, Z, \Delta, \Psi)|^{-1} = -2\pi\{\sum \delta_l^2 (\log \psi_l)^{-1} [1 + (\log \psi_l)^{-1} \log R_l] + \sum \delta_j \delta_l (\log \psi_j)^{-1} (\log \psi_l)^{-1} g_j(z_l)\}. \quad (7)$$

Relation (7) combined with Lemma 1 yields the required formula (6).

In the case of an arbitrary set $B$, the domains $B_l$ are exhausted by sequences of domains $B_l^s$ having the classical Green function:

$$\overline{B}_l^s \subset B_l^{s+1} \subset B_l, \quad s = 1, 2, \dots, \quad \bigcup_{s=1}^{\infty} B_l^s = B_l.$$

Let us introduce some notation:

$$B^s = \bigcup_{l=1}^{m} B_l^s, \quad g_l^s(z) = g_{B_l^s}(z, z_l), \quad R_l^s = r(B_l^s, z_l).$$

It is known that $g_l^s \to g_l$ as $s \to \infty$ uniformly inside $B_l \setminus z_l$ and, in particular, $R_l^s \to R_l$ as $s \to \infty$, $l = 1, \dots, m$. Let $g(s)$ and $\widetilde{E}(z_l, r)$ be the same as before, and let $g^s(z)$ be the function constructed from $g_l^s(z)$ similarly to $g(z)$. Standard arguments using the uniform convergence show that

$$\lim_{s\to\infty} \iint\limits_{B^s\setminus \bigcup_{l=1}^{m} \widetilde{E}(z_l,r)} |\nabla g^s|^2 \, dxdy = |\widetilde{C}(r; B, Z, \Delta, \Psi)|^{-1}.$$

On the other hand, we see as before that the integral under the limit sign is equal to

$$-\sum_{j=1}^{m} \int\limits_{\partial \widetilde{E}(z_j,r)} (\delta_j + \alpha_{js}(z)) \frac{\partial g^s}{\partial n} ds = -\sum_{j=1}^{m} \delta_j \int\limits_{\partial E(z_j,\rho)} \frac{\partial g^s}{\partial n} ds + \alpha_s$$

$$= -2\pi \left\{ \sum \delta_l^2 (\log \psi_l)^{-1} [1 + (\log \psi_l)^{-1} \log R_l^s] + \sum \delta_j \delta_l (\log \psi_j)^{-1} (\log \psi_l)^{-1} g_j^s(z_l) \right\} + \alpha_s,$$

where $\alpha_{js}(z)$ and $\alpha$ are some infinitesimals as $s \to \infty$. Passing to the limit as $s \to \infty$, we once more obtain relation (7) and, hence, formula (6). The theorem is proved.

## §2. An application to geometric function theory

First, note that if $B$ and $B'$ are open sets of the plane $\mathbf{C}_z$ and $B \subset B'$, then the monotonicity of the reduced module and Theorem 1 imply the Nehari inequality (see [6])

$$\sum \delta_l \delta_j h_l(z_l, z_j) \leq \sum \delta_l \delta_j h_l'(z_l, z_j) \tag{8}$$

for any real $\delta_l$, $l = 1, \ldots, m$. Here $h_l(z,\zeta) = g_{B_l}(z,\zeta) + \log|z-\zeta|$ and $h_l'(z,\zeta) = g_{B_l'}(z,\zeta) + \log|z-\zeta|$ are the regular parts of the Green functions of the connected components of $B$ and $B'$, respectively (see also [7, §3]). Setting $B' = \{z : |z| < R\}$ in (8) and passing to the limit as $R \to \infty$, we obtain the well-known inequality

$$\prod_{l=1}^{m} r(B_l, z_l)^{\delta_l^2} \leq \prod_{1\leq k<l\leq m} |z_k - z_l|^{-2\delta_k\delta_l},$$

where $\sum \delta_l = 0$ and $B_l$, $l = 1, \ldots, m$, are pairwise nonoverlapping domains [5, p. 551].

Now let us proceed to new results obtained by using dissymmetrization and the separating transformation. We will stick to the terminology used in [4]. The only distinction is that in our case, generally speaking, condensers consist of a large number of plates. However, the definitions are either analogous, or exactly repeated.

Let $\Phi_n$ be the group of symmetries of $\overline{\mathbf{C}}_z$ which consists of superpositions of reflections with respect to the lines $\{z = t\exp(\pi ik/n) : -\infty < t < +\infty\}$, $k = 1, \ldots, n$. A set $A$ is called *$\Phi_n$-symmetric* if $\varphi(A) = A$ for all isometries $\varphi \in \Phi_n$. We call a collection of points *$\Phi_n$-symmetric* if the set constituted by these points is symmetric. The orbit of the point $z$ with respect to the group $\Phi_n$ is the set $\{\varphi(z) : \varphi \in \Phi_n\}$.

**Theorem 2.** *Let $B^*$ be an open $\Phi_n$-symmetric set on the complex sphere $\overline{\mathbf{C}}_z$, let $\{z_l^*\}$ be a $\Phi_n$-symmetric collection of $m$ distinct points of the set $B^*$, and let $\{\delta_l\}$ be a collection of $m$ real numbers satisfying the following condition: $\delta_k = \delta_l$ if the points $z_k^*$ and $z_l^*$ belong to one and the same orbit. Let the set $B$ and the points $\{z_l\}$ be the results of a certain dissymmetrization of $B^*$ and $\{z_l^*\}$, respectively. Then we have the inequality*

$$\sum \delta_l^2 \log r(B_l^*, z_l^*) + \sum \delta_k\delta_l g_{B_l^*}(z_k^*, z_l^*) \leq \sum \delta_l^2 \log r(B_l, z_l) + \sum \delta_k\delta_l g_{B_l}(z_k, z_l).$$

*Proof.* Consider the condenser $C(r; B^*, \{z_l^*\}, \{\delta_l\})$. It is easy to see that for sufficiently small $r > 0$ the result of dissymmetrization of this condenser coincides with the condenser $C(r; B, \{z_l\}, \{\delta_l\})$. It follows from the dissymmetrization principle (similarly to Theorem 1.10 of [4]) that

$$|C(r; B^*, \{z_l^*\}, \{\delta_l\})| \leq |C(r; B, \{z_l\}, \{\delta_l\})|.$$

Passing to the limit as $r \to 0$, we obtain

$$M(B^*, \{z_l^*\}, \{\delta_l\}) \leq M(B, \{z_l\}, \{\delta_l\}).$$

It remains to use the formula for the reduced module of Theorem 1.

Let us mention some corollaries of Theorem 2. Let $a = (a_1, \ldots, a_n)$ be a collection of arbitrary distinct points on the circle $|z| = 1$, and let $a^* = (1, \exp(2\pi i/n), \ldots, \exp(2\pi i(n-1)/n))$. Let $K$ be a closed set on the segment $(0, 1]$, and let $g_a(z, \zeta)$ be the Green function of the domain

$$B_a = \{z : |z| < 1\} \setminus \bigcup_{k=1}^{n} \{z = a_k t : t \in K\}.$$

It is known that there exists a dissymmetrization taking the set $B_{a^*}$ to the set $B_a$ [4, p. 34]. Let us fix $r$, $0 < r < 1$, and consider the $\Phi_n$-symmetric collection of points $z_l^* = r\exp(i(\pi/(nN) + 2\pi(l-1)/(nN)))$.

$l = 1, \dots, nN$, and the point $z = 0$. The corresponding collection of real numbers is constructed by the numbers $\delta_l = 1/(nN)$, $l = 1, \dots, nN$, and an arbitrary real number $\delta$. Theorem 2 yields the inequality

$$\frac{1}{(nN)^2} \sum \log r(B_{a^*}, z_l^*) + \delta^2 \log r(B_{a^*}, 0) + 2\delta \frac{1}{nN} \sum g_{a^*}(0, z_l^*) + \frac{1}{(nN)^2} \sum g_{a^*}(z_k^*, z_l^*)$$

$$\leq \frac{1}{(nN)^2} \sum \log r(B_a, z_l) + \delta^2 \log r(B_a, 0) + 2\delta \frac{1}{nN} \sum g_a(0, z_l) + \frac{1}{(nN)^2} \sum g_a(z_k, z_l).$$

Passing to the limit as $N \to \infty$, we obtain the inequality

$$\delta^2 \log r(B_{a^*}, 0) + \frac{\delta}{\pi} \int_0^{2\pi} g_{a^*}(0, re^{i\theta}) d\theta + \frac{1}{(2\pi)^2} \int_0^{2\pi}\int_0^{2\pi} g_{a^*}(re^{i\varphi}, re^{i\theta})\, d\varphi d\theta$$

$$\leq \delta^2 \log r(B_a, 0) + \frac{\delta}{\pi} \int_0^{2\pi} g_a(0, re^{i\theta}) d\theta + \frac{1}{(2\pi)^2} \int_0^{2\pi}\int_0^{2\pi} g_a(re^{i\varphi}, re^{i\theta})\, d\varphi d\theta.$$

Setting here $\delta = 0$, we get

$$\int_0^{2\pi}\int_0^{2\pi} g_{a^*}(re^{i\varphi}, re^{i\theta})\, d\varphi d\theta \leq \int_0^{2\pi}\int_0^{2\pi} g_a(re^{i\varphi}, re^{i\theta})\, d\varphi d\theta.$$

The case $\delta \to \infty$, or $r \to 1$, gives the known inequality

$$r(B_{a^*}, 0) \leq r(B_a, 0).$$

Note that according to the Haliste conjecture concerning Gonchar's well-known problem [8], we have

$$\int_0^{2\pi} g_{a^*}(0, re^{i\theta}) d\theta \leq \int_0^{2\pi} g_a(0, re^{i\theta}) d\theta,$$

and up to now this inequality has been neither proved nor refuted.

Now consider arbitrary points $z_k$, $k = 1, \dots, n$, lying on the circle $|z| = r$, and the points $\zeta_k = z_k R/r$, $k = 1, \dots, n$, $0 < r$, $R < \infty$. Let $z_k^*$ and $\zeta_k^* = z_k^* R/r$, $k = 1, \dots, n$, constitute a $\Phi_n$-symmetric collection of points. It is known that there exists a dissymetrization taking the points $z_k^*$, $\zeta_k^*$ to the points $z_k$, $\zeta_k$ (see [4]). Applying Theorem 2 to these points and to the domain $B_\rho^* = B_\rho = \{z : |z| > \rho\}$ and then passing to the limit as $\rho \to 0$, we obtain the inequality

$$\frac{\prod |z_k^* - \zeta_l^*|^{2\alpha\beta}}{\prod |z_k^* - z_l^*|^{\alpha^2} |\zeta_k^* - \zeta_l^*|^{\beta^2}} \leq \frac{\prod |z_k - \zeta_l|^{2\alpha\beta}}{\prod |z_k - z_l|^{\alpha^2} |\zeta_k - \zeta_l|^{\beta^2}},$$

which holds for all real $\alpha$ and $\beta$. For $\alpha = 1$ and $\beta = 0$, we have the classical inequality (see [9])

$$\prod |z_k - z_l| \leq \prod |z_k^* - z_l^*|.$$

It is clear that a nontrivial lower estimate of the left-hand side of this inequality is possible only under additional restrictions on the points $z_k$. In this connection, the following result is of interest.

**Theorem 3.** *Let $z_{kl} \neq 0$ be distinct points lying on radial rays so that $\arg z_{kl} = 2\pi k/n$, $l = 1, \dots, m_k$, $k = 1, \dots, n$, and let the real numbers $\delta_{kl}$, $l = 1, \dots, m_k$, $k = 1, \dots, n$, satisfy the conditions*

$$\nu^2 \left(\sum \delta_{kl}\right)^2 = \frac{1}{n} \sum_{k=1}^{n} \nu_k^2 \left(\sum_{l=1}^{m_k} \delta_{kl}\right)^2, \quad \nu = \frac{1}{n^2} \sum \nu_k, \tag{9}$$

*where*

$$\nu = \left(\sum \delta_{kl}^2\right)^{-1}, \quad \nu_k = \left(\sum_{l=1}^{m_k} \delta_{kl}^2\right)^{-1}, \quad k = 1, \dots, n.$$

*Then we have the inequality*

$$\prod |z_{ki} - z_{lj}|^{n^2 \nu^2 \delta_{ki} \delta_{lj}} \geq \prod (n|z_{kl}|^{n-1})^{\nu^2 \delta_{kl}^2} \prod \big| |z_{ki}|^n - |z_{kj}|^n \big|^{\nu_k^2 \delta_{ki} \delta_{kj}}. \tag{10}$$

*Proof.* Let us introduce the following notation:

$$B_\rho = \{z : |z| < \rho\}, \quad Z = \{z_{kl}\}, \quad \Delta = \{\delta_{kl}\},$$

$$w = p_k^+(z) \equiv -iz^n, \quad \frac{2\pi k}{n} < \arg z < \frac{2\pi k}{n} + \frac{\pi}{n}, \quad k = 1, \dots, n,$$

$$w = p_k^-(z) \equiv -iz^n, \quad -\frac{\pi}{n} + \frac{2\pi k}{n} < \arg z < \frac{2\pi k}{n}, \quad k = 1, \dots, n.$$

For sufficiently large $\rho$ and small $r > 0$, let $\{C_k^+, C_k^-\}_{k=1}^n$ be the result of the separating transformation of the condenser $C(r; B_\rho, Z, \Delta)$ with respect to the family of functions $\{p_k^+(z), p_k^-(z)\}_{k=1}^n$ (see [4, p. 27]). It is easy to see that the plates of the condensers $C_k^\pm$ satisfy conditions (3), so that

$$C_k^\pm = \widetilde{C}\left(r; B_{\rho^n}, \{p_k^\pm(z_{kl})\}_{l=1}^{m_l}, \{\delta_{kl}\}_{l=1}^{m_k} \{n|z_{kl}|^{n-1} r\}_{l=1}^{m_k}\right).$$

As in Theorem 1.8 of [4], we have

$$|C(r; B_\rho, Z, \Delta)|^{-1} \geq \frac{1}{2} \sum_{k=1}^{n} \left(|C_k^+|^{-1} + |C_k^-|^{-1}\right).$$

Taking into account that $|C_k^+| = |C_k^-|$, we obtain

$$|C(r; B_\rho, Z, \Delta)| \leq \frac{1}{n^2} \sum_{k=1}^{n} |C_k^+|.$$

In view of (9), this implies

$$|C(r; B_\rho, Z, \Delta)| + \frac{\nu}{2\pi} \log r \leq \frac{1}{n^2} \left\{ \sum_{k=1}^{n} \left[ |C_k^+| + \frac{\nu_k}{2\pi} \log r \right] \right\}.$$

By Theorem 1, we conclude that

$$\frac{\nu^2}{2\pi} \left\{ \sum \delta_{kl}^2 \log \frac{\rho^2 - |z_{kl}|^2}{\rho} + \sum \delta_{ki} \delta_{lj} \log \left| \frac{\rho^2 - \bar{z}_{ki} z_{lj}}{\rho(z_{ki} - z_{lj})} \right| \right\} l$$

$$\leq \frac{1}{n^2} \sum_{k=1}^{n} \frac{\nu_k^2}{2\pi} \left\{ \sum_{l=1}^{m_k} \delta_{kl}^2 \log \frac{\rho^{2n} - |z_{kl}|^{2n}}{\rho^n n |z_{kl}^{n-1}} + \sum_{ij} \delta_{ki} \delta_{kj} \log \left| \frac{\rho^{2n} - |z_{ki} z_{kj}|^n}{\rho^n (|z_{ki}|^n - |z_{kj}|^n)} \right| \right\}.$$

It remains to use relation (9) and pass to the limit as $\rho \to 0$. The theorem is proved.

If each one of the rays $\arg z = 2\pi k/n$ contains exactly one of the points $z_k$, $k = 1, \dots, n$, then inequality (10) yields

$$\prod |z_k - z_l| \geq n^n \prod |z_k|^{n-1} \tag{11}$$

($z_{k1} = z_k$, $m_k = 1$, $\delta_{k1} = 1$, $\nu_k = 1$, $\nu = 1/n$). Here, the equality is attained for $\Phi_n$-symmetric points $z_k$. In the case where each of the rays contains two of the points, we have, e.g., the inequality

$$\frac{n^{2n} \prod |z_{k1} - z_{l2}|}{\prod |z_{k1} - z_{l1}|\, |z_{k2} - z_{l2}|} \leq \prod \frac{\big| |z_{k1}|^n - |z_{k2}|^n \big|^2}{|z_{k1} z_{k2}|^{n-1}}$$

($\delta_{k1} = -\delta_{k2} = 1$, $k = 1, \dots, n$). In addition to these inequalities, let us indicate a way of strengthening them. Here we restrict ourselves to the simplest example, i.e., to inequality (11). Let

$$z_k \neq 0, \quad \arg z_k = 2\pi k/n, \quad k = 1, \dots, n, \quad Z = \{z_k\}, \quad \Delta = \{1, \dots, 1\},$$

$$w = p_k(z) \equiv -i z^{n/2}, \quad \frac{2\pi k}{n} < \arg z < \frac{2\pi(k+1)}{n}, \quad p_k(z_k) = -i|z_k|^{n/2}, \quad k = 1, \dots, n.$$

Let $\{C_k\}$ be the result of the separating transformation of the condenser $C(r; B_\rho, Z, \Delta)$ with respect to the family of functions $\{p_k(z)\}$. By Theorem 1.8 of [4], we obtain

$$|C(r; B_\rho, Z, \Delta)| + \frac{1}{2\pi n} \log r \leq \frac{2}{n^2} \sum \left[ |C_k| + \frac{1}{4\pi} \log r \right].$$

As before, the asymptotic formula of Theorem 1 leads to the inequality

$$\prod |z_k - z_l| \geq \left(\frac{n}{2}\right)^n \prod |z_k|^{n/2-1} \left(|z_k|^{n/2} + |z_{k+1}|^{n/2}\right) \quad (z_{n+1} = z_1).$$

Problems of extremal decomposition that involve modules and reduced modules of domains of the form $M(B, z_0)$ are well known in function theory (see [10]). We continue the study of similar problems for generalized reduced modules with respect to collections of points [2].

**Theorem 4.** *Let $D_k$, $k = 0, 1, \dots, n$, be domains of the complex sphere $\overline{\mathbf{C}}_z$ such that $D_0$ contains the point $z = 0$ and does not intersect the domains $D_k$, where $D_k$ contains a point $z_k$ with $\arg z_k = 2\pi k/n$, $k = 1, \dots, n$. Consider the connected component of $\bigcup_{k=1}^n D_k$ containing the point $z_k$. Let $g_k(z)$ be its Green function with pole at this point, and let $g_k(z) = 0$ outside of this connected component. Then*

$$\exp\Big(2 \sum_{1 \leq k < l \leq n} g_k(z_l)\Big) (r(D_0, 0))^{n^2} \prod_{k=1}^n r(D_k, z_k) \leq n^{-n} \prod_{k=1}^n |z_k|^{n+1}. \tag{12}$$

*The equality here is attained for the domains $D_0^*$ and $D_k^*$ bounded by the curves $z = \sqrt[n]{R^n/2 + it}$, $-\infty < t < +\infty$, and the points $z_k^* = R\exp(2\pi i k/n)$, $k = 1, \dots, n$, where $R$ is an arbitrary positive number.*

*Proof.* Put

$$B_\rho = \{z : |z| < \rho\}, \quad \rho > \max_{1 \leq k \leq n} |z_k|, \quad Z = \{0, z_1, \dots, z_n\}, \quad \Delta = \{-n, 1, \dots, 1\}.$$

Let $D_{k\rho}$ be the connected component of $D_k \cap B_\rho$ containing the point $z_k$ ($z_0 = 0$). Consider the connected component of $\bigcup_{k=1}^n D_k \cap B_\rho$ containing the point $z_k$, and let $g_{k\rho}(z)$ be its Green function with pole at this point. By Theorem 1 and in view of the monotonicity of the reduced module, we have

$$\frac{n^2 \log r(D_{0\rho}, 0) + \sum_{k=1}^n \log r(D_{k\rho}, z_k) + 2 \sum_{1 \leq k < l \leq n} g_{k\rho}(z_l)}{2\pi (n(n+1))^2} \leq M\Big(\bigcup_{k=0}^N D_{k\rho}, Z, \Delta\Big) \leq M(B_\rho, Z, \Delta).$$

Direct calculations yield

$$M(B_\rho, Z, \Delta) = \frac{1}{2\pi(n(n+1))^2}\left\{n^2\log\rho + \sum\log\frac{\rho^2-|z_k|^2}{\rho} + \sum\log\left|\frac{\rho^2-\bar{z}_k z_l}{\rho(z_k-z_l)}\right| - 2n\sum\log\left|\frac{\rho}{z_k}\right|\right\}.$$

If we sum up the above relations and pass to the limit as $\rho \to \infty$, then the obtained inequality combined with (11) yields the required result (12). The equality case is verified directly. The theorem is proved.

Note that $D_k^*$, $k = 1,\dots,n$, are simply connected pairwise disjoint domains. Thus, inequality (12) generalizes the corresponding inequality in the problem of extremal decomposition with free poles in the classical setting [10] to the case of partly nonoverlapping domains (cf. [4, §9]). An effective method of solving problems of extremal decomposition consists of using a separating transformation of domains [4]. A separating transformation of arbitrary sets and the behavior of generalized reduced modules under such a transformation is the subject of a special work. For this reason, here we restrict ourselves to only two canonical cases of extremal decomposition, where the generalized reduced module reduces to the usual one, and hence the extremal problem reduces to the familiar problem of extremal decomposition. For this purpose, we need the following lemma.

**Lemma 2.** *Assume that the domain $D$ has the classical Green function $g(z,\zeta)$, and let $z_1$ and $z_2$ be two points of this domain. Then there exist disjoint domains $D_1$ and $D_2$ such that $D_1 \cup D_2 \subset D$ and*

$$\log(r(D_1,z_1)r(D_2,z_2)) = \log(r(D,z_1)r(D,z_2)) - 2g(z_1,z_2).$$

*Proof.* The function $u(z) = g(z,z_1) - g(z,z_2)$ is defined on the domain $D \setminus \{z_1,z_2\}$ and is harmonic in it; moreover, $u(z) \to +\infty$ as $z \to z_1$ and $u(z) \to -\infty$ as $z \to z_2$. It follows that the sets $D_1 = \{z : u(z) > 0\}$ and $D_2 = \{z : u(z) < 0\}$ are domains. The Green function of the domain $D_1$ with pole at the point $z_1$ coincides with the function $u(z)$. Therefore,

$$\log r(D_1,z_1) = \lim_{z\to z_1}\left(u(z) + \log|z-z_1|\right) = \log(D,z_1) - g(z_1,z_2).$$

In a similar way, $-u(z)$ coincides with the Green function of the domain $D_2$ with pole at $z_2$ and

$$\log r(D_2,z_2) = \log r(D,z_2) - g(z_2,z_1).$$

It remains to sum up the relations obtained. The lemma is proved.

**Theorem 5.** *For any points $z_k$, $\zeta_k$ satisfying the condition*

$$\arg z_k = \arg\zeta_k, \quad |z_k| = r, \quad |\zeta_k| = R, \quad k = 1,\dots,n, \quad 0 < r < R < \infty,$$

*and any pairwise nonoverlapping domains $D_k$ such that $z_k$, $\zeta_k \in D_k \subset \overline{\mathbb{C}}_z$, $k = 1,\dots,n$, we have the inequality*

$$\prod_{k=1}^{n} r(D_k,z_k)r(D_k,\zeta_k)\exp(-2g_k(z_k,\zeta_k)) \le \left[\frac{4\sqrt{rR}(R^{n/2}-r^{n/2})}{n(R^{n/2}+r^{n/2})}\right]^{2n},$$

*where $g_k(z,\zeta)$ is the Green function of the domain $D_k$. The equality is attained for the domains*

$$D_k^* = \left\{z : \frac{2\pi k}{n} < \arg z < \frac{\pi(k+1)}{n}\right\}$$

*and for the points*

$$z_k^* = r\exp\left(i\left(\frac{\pi}{n} + \frac{2\pi k}{n}\right)\right), \quad \zeta_k^* = R\exp\left(i\left(\frac{\pi}{n} + \frac{2\pi k}{n}\right)\right), \quad k = 1,\dots,n.$$

*Proof.* We can assume that the domain $D_k$ possesses a Green function. By Lemma 2, for each domain $D_k$ there exist two nonoverlapping domains $D_{k1}$ and $D_{k2}$ such that $D_{k1} \cup D_{k2} \subset D_k$ and

$$r(D_k, z_k) r(D_k, \zeta_k) \exp(-2g_k(z_k, \zeta_k)) = r(D_{k1}, z_k) r(D_{k2}, \zeta_k), \quad k = 1, \dots, n.$$

Hence, the proof of the theorem reduces to obtaining an upper estimate of the product

$$\prod_{k=1}^{n} r(D_{k1}, z_k) r(D_{k2}, \zeta_k), \tag{13}$$

where $D_{k1}$, $D_{k2}$, $k = 1, \dots, n$, are pairwise nonoverlapping domains. The easiest way to obtain such an estimate consists of the successive application of a separating transformation and Theorem 2.15 of [4]. On the other hand, Emel'yanov has already obtained the indicated estimate for simply connected domains (see [11, p. 96]). His method of proof is valid in this case with adjustments known to specialists, and also for multiply connected domains. The key feature here is the specific character of quadratic differential (12) used in [11]. Repeating Emel'yanov's proof, we obtain

$$\prod_{k=1}^{n} r(D_{k1}, z_k) r(D_{k2}, \zeta_k) \le \left[ \frac{4\sqrt{rR}(R^{n/2} - r^{n/2})}{n(R^{n/2} + r^{n/2})} \right]^{2n}.$$

The verification of the equality case is an easy exercise. The theorem is proved.

**Theorem 6.** *For any points $z_k$, $\zeta_k$ lying on the circle $|z| = r$ and any pairwise nonoverlapping domains $D_k$ such that $z_k, \zeta_k \in D_k \subset \overline{\mathbf{C}}_z$, $k = 1, \dots, n$, we have the inequality*

$$\prod_{k=1}^{n} r(D_k, z_k) r(D_k, \zeta_k) \exp(-2g_k(z_k, \zeta_k)) \le (2r/n)^{2n}.$$

*The equality is attained for the domains*

$$D_k^* = \left\{ z : \frac{2\pi k}{n} < \arg z < \frac{2\pi(k+1)}{n} \right\}$$

*and the points*

$$z_k^* = r \exp\left( i\left( \frac{\pi}{2n} + \frac{2\pi k}{n} \right) \right), \quad \zeta_k^* = r \exp\left( i\left( \frac{3\pi}{2n} + \frac{2\pi k}{n} \right) \right), \quad k = 1, \dots, n.$$

*Proof.* As in the previous proof, we see that it is suffices to give an upper estimate for product (13), where this time the points $z_k$ and $\zeta_k$ lie on one circle. By Theorem 2.17 of [4], this product does not exceed $(2r/n)^{2n}$. The equality case is verified by direct calculation of the product on the left for the indicated domains $D_k^*$ and points $z_k^*$, $\zeta_k^*$, $k = 1, \dots, n$. The theorem is proved.

The above results are directly connected with properties of regular and univalent functions. Let $B \subset \overline{\mathbf{C}}_z$ and $D \subset \overline{\mathbf{C}}_w$ be open sets such that their complements $\overline{\mathbf{C}}_z \setminus B$ and $\overline{\mathbf{C}}_w \setminus D$ have positive harmonic measure. Let the function $w = f(z)$ conformally and univalently map the set $B$ onto $D$. For any distinct points $Z = \{z_l\}$ of the set $B$ and any nonzero numbers $\Delta = \{\delta_l\}$, the reduced module $M(B, Z, \Delta)$ is defined. Define $W = \{f(z_l)\}$. By formula (6), we have

$$M(D, W, \Delta) = M(B, Z, \Delta) + \frac{\nu^2}{2\pi} \sum \delta_l^2 \log |f'(z_l)|. \tag{14}$$

Together with estimates of the modules of the condensers $C(r; B, Z, \Delta)$ and $C(r; D, W, \Delta)$, relation (14) gives applications of relation (6) to the theory of univalent functions.

**Theorem 7.** *Let $w = f(z)$ be a function which conformally and univalently maps the disk $|z| < 1$ onto a bounded domain $D$ such that the transfinite diameter of its closure is less than or equal to one. Assume that the points $w_0$ and $z_k$, $k = 1, \dots, n$, satisfy the following conditions: $f(z_k) \neq w_0$ and*

$$\arg \frac{f(z_1) - w_0}{f(z_k) - w_0} = \frac{2\pi(k-1)}{n}, \quad k = 1, \dots, n.$$

*Then*

$$n^n \prod |f'(z_k)|\, |f(z_k) - w_0|^{n-1} \le \frac{\prod |z_k - z_l|}{\prod |1 - z_k \bar{z}_l|}.$$

*Proof.* Let us introduce the following notation:

$$B = \{z : |z| < 1\}, \quad Z = \{z_k\}, \quad \Delta = \{\delta_k\}, \quad \delta_k = 1, \quad k = 1, \dots, n; \quad F(z) = \frac{1}{f(z) - w_0}; \quad W = \{F(z_k)\}.$$

By Theorem 1 and relation (14), we subsequently obtain

$$\frac{1}{2\pi n^2}\left\{\sum \log(1 - |z_k|^2) + \sum \log\left|\frac{1 - z_k \bar{z}_l}{z_k - z_l}\right|\right\} = M(B, Z, \Delta)$$

$$= -\frac{1}{2\pi n^2} \sum \log |f'(z_k)|\, |f(z_k) - w_0|^{-2} + M(F(B), W, \Delta).$$

Let $D_0$ be that connected component of the image of the set $\overline{\mathbb{C}}_w \setminus \overline{D}$ under the mapping $\zeta = 1/(w - w_0)$ which contains the origin. By the condition of the theorem, we have $r(D_0, 0) \ge 1$. Let $g(\zeta, \zeta')$ denote the Green function of the domain $F(B)$. Once again using Theorem 1 and inequality (12), we obtain

$$\begin{aligned} M(F(B), W, \Delta) &= \frac{1}{2\pi n^2}\left\{\sum \log r(F(B), F(z_k)) + \sum g(F(z_k), F(z_l))\right\} \\ &\le \frac{1}{2\pi n^2}\left\{n^2 \log r(D_0, 0) + \sum \log r(F(B), F(z_k)) + \sum g(F(z_k), F(z_l))\right\} \\ &\le \frac{1}{2\pi n^2} \log\left\{n^{-n} \prod |F(z_k)|^{n+1}\right\} = \frac{1}{2\pi n^2} \log\left\{n^{-n} \prod |f(z_k) - w_0|^{-n-1}\right\}. \end{aligned}$$

Summing up these relations, we get the required inequality. The theorem is proved.

**Remark.** This proof involves the notion of reduced module only for reasons of a methodical character. Our aim is to show the technique of applying this notion. The substantive part of the proof consists of applying inequality (12), which can be considered as an estimate of the reduced module. In the following proof, we will use the result (Theorem 5) and the behavior of the inner radius and the Green function under a conformal mapping.

Let

$$S_f(z) = \frac{f'''(z)}{f'(z)} - \frac{3}{2}\frac{(f''(z))}{(f'(z))^2}$$

be the Schwarz derivative of the function $f(z)$ at the point $z$.

**Theorem 8.** *For any functions $w = f_k(z)$ which are meromorphic and univalent in the disk $|z| < 1$ and map this disk onto mutually nonoverlapping domains so that $|f_k(0)| = 1$, $k = 1, \dots, n$, we have the sharp inequality*

$$\sum_{k=1}^{n} \left\{ |f_k'(0)|^{-2} - \frac{1}{6} \operatorname{Re} \frac{f_k^2(0) S_{f_k}(0)}{(f_k'(0))^2} \right\} \ge \frac{n(n^2 + 2)}{24}.$$

The equality here is attained only for the functions

$$f_k^*(z) = \exp\left(\frac{2\pi i(k-1)}{n}\right) \frac{[(1+z)/(1-z)]^2}{n},$$

where we take the branch of the root preserving the unity, $k = 1, \dots, n$.

*Proof.* Set

$$w_k = f_k(0), \quad D_k = f_k(\{z : |z| < 1\}), \quad w_k^* = f_k^*(0), \quad D_k^* = f_k^*(\{z : |z| < 1\}), \quad k = 1, \dots, n.$$

As in [3], it is not difficult to check that

$$|f_k'(0)|^{-2} - \frac{1}{6} \operatorname{Re} \frac{w_k^2 S_{f_k}(0)}{(f_k'(0))^2} = K(D_k, w_k), \tag{15}$$

where

$$K(D_k, w_k) = \lim_{\rho \to 0} \frac{-1}{4\rho^2} \{ \log r(D_k, (1-\rho)w_k) + \log r(D_k, (1+\rho)w_k) - 2g_{D_k}((1-\rho)w_k, (1+\rho)w_k) - 2\log(2\rho) \}.$$

By Theorem 5, we have

$$\sum_{k=1}^{n} K(D_k, w_k) \geq \sum_{k=1}^{n} K(D_k^*, w_k^*).$$

It remains to use relation (15). The theorem is proved.

This research was partially supported by the Russian Foundation for Basic Research, grant 96-01-00007.

Translated by N. Yu. Netsvetaev.